\documentclass{article}

\usepackage{amsfonts}
\usepackage{amssymb}
\usepackage{enumerate}
\usepackage{amsmath}
\usepackage[T1]{fontenc}

\newtheorem{ttt}{Theorem}[section]
\newtheorem{llll}[ttt]{Lemma}
\newtheorem{ccc}[ttt]{Claim}
\newtheorem{eee}[ttt]{Example}
\newtheorem{fff}[ttt]{Fact}
\newtheorem{rrr}[ttt]{Remark}
\newtheorem{sss}[ttt]{Statement}
\newtheorem{ddd}[ttt]{Definition}
\newtheorem{qqq}[ttt]{Question}
\newtheorem{cccc}[ttt]{Corollary}
\newtheorem{nnn}[ttt]{Notation}

\newcommand{\bt}{\begin{ttt}}
\newcommand{\bl}{\begin{llll}}
\newcommand{\bc}{\begin{ccc}}
\newcommand{\bex}{\begin{eee}}
\newcommand{\bfa}{\begin{fff}}
\newcommand{\br}{\begin{rrr}\upshape}
\newcommand{\bs}{\begin{sss}}
\newcommand{\bd}{\begin{ddd}\upshape}
\newcommand{\bq}{\begin{qqq}}
\newcommand{\bnn}{\begin{nnn}}
\newcommand{\bcor}{\begin{cccc}}

\newcommand{\bp}{\noindent\textbf{Proof. }}

\newcommand{\et}{\end{ttt}}
\newcommand{\el}{\end{llll}}
\newcommand{\ec}{\end{ccc}}
\newcommand{\eex}{\end{eee}}
\newcommand{\efa}{\end{fff}}
\newcommand{\er}{\end{rrr}}
\newcommand{\es}{\end{sss}}
\newcommand{\ed}{\end{ddd}}
\newcommand{\eq}{\end{qqq}}
\newcommand{\ecor}{\end{cccc}}
\newcommand{\enn}{\end{nnn}}

\newcommand{\ep}{\hspace{\stretch{1}}$\square$\medskip}

\newcommand{\lab}[1]{\label{#1}}

\newcommand{\NN}{\mathbb{N}}

\newcommand{\QQ}{\mathbb{Q}}
\newcommand{\RR}{\mathbb{R}}

\newcommand{\al}{\alpha}

\newcommand{\de}{\delta}
\newcommand{\e}{\varepsilon} 
\newcommand{\ph}{\varphi} 
\newcommand{\om}{\omega}
\newcommand{\si}{\sigma}

\newcommand{\sm}{\setminus}
\newcommand{\eset}{\emptyset}
\newcommand{\beq}{\begin{equation}}
\newcommand{\eeq}{\end{equation}}
\def\su{\subset}

\newcommand\rest{\mathord{\upharpoonright}}

\title{On a converse to Banach's Fixed Point Theorem}

\author{M\'arton Elekes\thanks{Partially supported by Hungarian Scientific
Foundation grants no.~49786, 37758 and F 43620.}}

\begin{document}

\maketitle 

\begin{abstract}
We say that a metric space $(X,d)$ possesses the \emph{Banach Fixed Point
Property (BFPP)} if every contraction $f:X\to X$ has a fixed point. The Banach
Fixed Point Theorem says that every complete metric space has the
BFPP. However, E.~Behrends pointed out \cite{Be1} that the converse
implication does not hold; that is, the BFPP does not imply completeness, in
particular, there is a non-closed subset of $\RR^2$ possessing the BFPP. He
also asked \cite{Be2} if there is even an open example in $\RR^n$, and whether
there is a `nice' example in $\RR$. In this note we answer the first question
in the negative, the second one in the affirmative, and determine the simplest
such examples in the sense of descriptive set theoretic complexity.

Specifically, first we prove that if $X\su\RR^n$ is open or $X\su\RR$ is
simultaneously $F_\si$ and $G_\de$ and $X$ has the BFPP then $X$ is
closed. Then we show that these results are optimal, as we give an $F_\si$ and
also a $G_\de$ non-closed example in $\RR$ with the BFPP.

We also show that a nonmeasurable set can have the BFPP. Our non-$G_\de$
examples provide metric spaces with the BFPP that cannot be remetrised by any
compatible complete metric. All examples are in addition bounded.
\end{abstract}

\insert\footins{\footnotesize{MSC codes: Primary 54H25, 47H10, 55M20, 03E15,
54H05 Secondary 26A16}}
\insert\footins{\footnotesize{Key Words: Banach, contraction, complete,
closed, Borel, sigma, delta, ambiguous, descriptive, transfinite, Lipschitz,
typical compact}}

\section{Introduction}

Converses to the Banach Fixed Point Theorem have a very long history. The
earliest such result seems to be that of Bessaga \cite{Bes}, but see also
\cite{Ba}, \cite{Ci}, \cite{Ja}, \cite{Ja1}, \cite{Ja2}, \cite{Ki}, \cite{Me},
\cite{MP}, \cite{Op} and \cite{Ru}. There are also numerous result of this
kind in linear spaces as well.

The version we consider in this note is the following.

\bd
We say that a metric space $(X,d)$ possesses the \emph{Banach Fixed Point
Property (BFPP)} if every contraction $f:X\to X$ has a fixed point.
\ed

Note that the empty set does not possess the BFPP as the empty function is a
contraction with no fixed point, so this would cause no
problem, but for the sake of simplicity we simply assume that all sets and
metric spaces considered are nonempty.

At the Problem Session of the 34th Winter School in Abstract Analysis
E.~Behrends presented the following example, which he referred to as
`folklore'.

\bt\lab{sin}
Let $X=graph \left(sin(1/x)\rest_{(0,1]} \right)$. Then $X\su\RR^2$ is a
non-closed set possessing the Banach Fixed Point Property.
\et

\bp
$X$ is clearly not closed. Let $f:X\to X$ be a contraction of Lipschitz
constant $q<1$.
For $H\su(0,1]$ define $X\rest_H=graph \left(sin(1/x)\rest_H \right)$. 
Choose $\e>0$ so that $diam \left( X\rest_{(0,\e)} \right) <
\frac{2}{q}$, then $diam \left( f\left( X\rest_{(0,\e)} \right) \right) <
2$. Hence $f\left( X\rest_{(0,\e)} \right)$ cannot contain both a local
minimum and a local maximum on the graph. 
But this set is clearly connected, which easily
implies that it is contained in at most two monotone parts of the
graph. Therefore there exists $\de_1>0$ such that $f\left( X\rest_{(0,\e)}
\right) \su X\rest_{[\de_1,1]}$. By compactness $f\left( X\rest_{[\e,1]}
\right) \su X\rest_{[\de_2,1]}$ for some $\de_2>0$, and hence setting 
$\de=min\{\de_1,\de_2\}$ gives $f(X) \su X\rest_{[\de,1]}$. But then the
Banach Fixed Point Theorem applied to $X\rest_{[\de,1]}$ provides a fixed point.
\ep

E.~Behrends asked the following two questions.

\bq\lab{q:1}\emph{(\cite{Be2})}
Is there an open non-closed subset of $\RR^n$ with the Banach Fixed Point
Property for some $n\in \NN$?
\eq

\bq\lab{q:2}\emph{(\cite{Be2})}
Is there a `simple' non-closed subset of $\RR$ with the Banach Fixed Point
Property?
\eq

\section{When Banach's Fixed Point Theorem implies completeness}

First we answer Question \ref{q:1}.

\bl\lab{segment}
Let $n\in\NN$ and $X\su\RR^n$ such that there exist $y,z\in\RR^n$ so
that $y\notin X$ but the nondegenerate segment $(y,z) \su X$. Then $X$ does
not possess the Banach Fixed Point Property.
\el

\bp
We can clearly assume $y=(0,\dots,0)$ and $z=(1,\dots,0)$. 
Then
\[
f(x) = \left(\frac{1}{2} arctan |x|,0,\dots,0\right)   \ \ \  \left( x\in
\RR^n\right) 
\]
is a contraction, since the absolute
value of vectors and $arctan$ are both Lipschitz functions of constant $1$. By our
assumptions $f(X)\su X$. As no contraction can have more than one fixed point,
and the origin is clearly a fixed point, we obtain that $f\rest_X$ has no
fixed point.
\ep

\bcor\lab{open}
For every $n\in\NN$ every open subset of $\RR^n$ possessing the Banach Fixed
Point Property coincides with $\RR^n$, hence it is closed.
\ecor

\bp
Let $U\su\RR^n$ be open but not closed, then there exists $z\in U$ and
$x\notin U$. Let $y$ be the closest point of $[x,z)\sm U$ to $z$.
\ep

Now we turn to Question \ref{q:2}, the case of $X\su\RR$. In this section we
show that there is no example that is simultaneously $F_\si$ and $G_\de$.
\bl
Let $X\su\RR$ such that $0\in \overline{X} \sm X$ and $0$ is a bilateral
accumulation point of $int(X^c)$. Then $X$ does not possess the Banach Fixed
Point Property.
\el

\bp
Let $\{(a_n,b_n)\}_{n\in\NN}$ be a sequence of intervals in $X^c \cap
(0,\infty)$ so that $b_{n+1} < a_n$ for every $n$ and $a_n,b_n\to 0$. Fix a
monotone decreasing sequence $z_n\in X$ such that $|z_n| <
\frac{b_n-a_n}{2}$. Now, for $x\in X$, $x>0$ let $n_x$ be
the minimal number for which $b_{n_x} < x$, and define $f(x) = z_{n_x}$. 
Define $f$ on $X\cap (-\infty,0)$ in a a similar manner.
We claim that $f$ is a contraction. 
First let $0<x<y$ be two points in $X$. If $n_x=n_y$ then $f(x) = f(y)$,
while 
if $n_x > n_y$ then $|f(x)-f(y)| < |z_{n_y}| < \frac{b_{n_y}-a_{n_y}}{2}
< \frac{|x-y|}{2}$, hence $f$ is a contraction on $X\cap
(0,\infty)$. Similarly, $f$ is a contraction on $X\cap (-\infty,0)$. Moreover,
$|f(x)| = |z_{n_x}| <
\frac{b_{n_x}-a_{n_x}}{2} < \frac{|x|}{2}$, which shows that for every
$x<0<y$ in 
$X$ we have $|f(x)-f(y)| < \frac{|x-y|}{2}$, hence $f$ is a contraction on
$X$.

Since $0 \notin X$, the above inequality $|f(x)| < \frac{|x|}{2}$ also shows
that $f$ has no fixed point. This finishes the proof.
\ep

A \emph{portion} of a set is a relatively open nonempty subset. A set that
is simultaneously $F_\si$ and $G_\de$ is called \emph{ambiguous} (or
$\Delta^0_2$ in descriptive set theory). 
A set $X$ is ambiguous iff for every nonempty closed set $F$
either $X$ or $X^c$ contains a portion of $F$ \cite{Ku}.

\bt\lab{amb}
Every simultaneously $F_\si$ and $G_\de$ subset of $\RR$ with the Banach Fixed
Point Property is closed.
\et

\bp
Suppose that $X\su\RR$ is a non-closed ambiguous set with the
BFPP. By applying a translation we can assume that $0\in \overline{X} \sm X$. 
By the previous lemma $0$ is not a bilateral accumulation point of
$int(X^c)$, so without loss of generality there exists $\e>0$ such that $X$ is
dense in $[0,\e]$. Let $I$ be an arbitrary closed nondegenerate subinterval of
$[0,\e]$. As $X$ is ambiguous, either $X$ or $X^c$ contains a portion
of $I$, but as $X$ is dense in $I$, the second alternative
cannot hold. Hence $X$ contains a subinterval of $I$, and as $I$ was
arbitrary, $int(X)$ is dense in $[0,\e]$. 

Set $F=[0,\e] \sm int(X)$. As $0\in F$, we have $F\neq\eset$, so either $X$ or
$X^c$ contains a portion 
of $F$, but the first alternative clearly cannot hold, so there exists an open
interval $J\su [0,\e]$ so that the nonempty set $F \cap J$ is disjoint from
$X$. Fix $f\in J\sm X$ and by the denseness of $int(X)$ also an $x\in J\cap
int(X)$. Let $y$ be the closest point to $x$ of $(int(X))^c$ between $x$ and
$f$. As $y\in F \cap J$, we obtain $y\notin X$, hence by Lemma \ref{segment} X
does not possess the BFPP. 
\ep

\section{When Banach's Fixed Point Theorem holds for strange sets}

In this section we give the examples of non-closed sets with the BFPP of
lowest possible Borel classes. For every $n\ge 2$ Theorem \ref{sin} clearly
provides an ambiguous example in $\RR^n$, 
Corollary \ref{open} shows that no open example
is possible, and obviously there is no closed example. In the language of
descriptive set theory, $\Delta^0_2$ is best possible, as there are no
$\Sigma^0_1$ and $\Pi^0_1$ examples. In $\RR$ Theorem \ref{amb} shows that
there is no ambiguous example, and this will be shown to be optimal when we
prove below that there are $F_\si$ and also $G_\de$ examples. That is,
$\Sigma^0_2$ and $\Pi^0_2$ are possible, but $\Delta^0_2$ is not.

The space of compact subsets of $\RR$ endowed with Hausdorff metric is a
complete metric space (see e.g.~\cite{Ke} for definitions and basic facts). 
We say that \emph{a typical compact set has a property} if the compact 
sets not having the property form a first category (in the sense of Baire) set
in the above space.

The following lemma is interesting in its own right. For simplicity we only
prove it in $\RR$, but it easily generalises to higher dimensions.

\bl\lab{cpt}
A contractive image of a typical compact $K\su\RR$ cannot contain a portion of
$K$.
\el

\bp
Recall that if each of a countable set of properties hold for a typical
compact set, then they also hold simultaneously, as first category sets are
closed under countable unions. Therefore it is enough to show that for a fixed
pair of rationals $p<q$, for a typical
compact set $K$ either $K\cap (p,q) = \eset$ or $K\cap (p,q)$ cannot be
covered by a contractive image of $K$. 
Similarly, it suffices to check that for a fixed
$r<1$ if $f$ is a contraction of ratio at most $r$ then either $K\cap (p,q) =
\eset$ or $K\cap (p,q) \not\subset f(K)$. 
As (in fact, in every dimension) every contraction can be extended to $\RR$
with the same Lipschitz constant \cite[2.10.43.]{Fe} we may assume that
$f:\RR\to\RR$.

Therefore it suffices to prove that for a fixed $r<1$ and for a fixed
pair of rationals $p<q$
\[
N = \{K\su \RR \textrm{ cpt }: \exists f:\RR\to\RR \textrm{
  contr.~of ratio} \le r, \ \ \eset \neq K\cap (p,q) \subset f(K)\}
\]
is a nowhere dense subset of the space of compact sets. 
Let $B(K_0,\e_0)$ be the open ball of center $K_0$ and radius $\e_0>0$. 
We need to find a ball inside this one that is disjoint from $N$. 
It is well known and easy to see that the finite sets form a dense subset of
our space, so we may assume that $K_0$ is finite; $K_0=\{x_1,\dots,x_n\}$.
%where these points are all distinct.

Suppose first that $K_0\cap [p,q] = \eset$. Define $\e_1 =
min\{dist(K_0,(p,q)), \e_0\} >0$. Then
for every $K\in B(K_0,\e_1)$ we have $K\cap (p,q) = \eset$, hence $B(K_0,\e_1)
\cap N = \eset$.

So we can assume that $K_0\cap [p,q] \neq \eset$, e.g.~$x_{i_0} \in [p,q]$. 
Let $(a,b)$ be a subinterval of $(p,q) \cap (x_{i_0}-\e_0,x_{i_0}+\e_0)$.
Choose an integer
\beq\lab{eq:k}
k > \frac{n+2}{1-r},
\eeq
and choose two arithmetic progressions $\{y_1,\dots,y_k\}$ and
$\{z_1,\dots,z_k\}$ in $(a,b)$, each of length $k$ and of some difference 
$d>0$ so that
\beq\lab{eq:dist}
dist(\{y_1,\dots,y_k\},\{z_1,\dots,z_k\}) \ge kd.
\eeq
Define
\[
K_1 = K_0 \cup \{y_1,\dots,y_k\} \cup \{z_1,\dots,z_k\},
\]
then $K_1\in B(K_0,\e_0)$. Choose 
\[
\e_1 = min\left\{dist\left(K_1, B(K_0, \e_0)^c\right), \frac{d}{4}\right\},
\]
then clearly $B(K_1,\e_1) \su B(K_0,\e_0)$. It is also easy to see that
the intervals $Y_i = (y_j-\e_1,y_j+\e_1)$, $Z_i = (z_j-\e_1,z_j+\e_1)$ for
$1\le j\le k$ are all disjoint. Also put $X_i = (x_i-\e_1, x_i+\e_1)$ for
every $1\le i\le n$.

Now we claim that $B(K_1,\e_1) \cap N = \eset$,
which will finish the proof. Let $K\in B(K_1,\e_1)$ be arbitrary. 
Clearly $K \su \bigcup_{i=1}^n X_i \cup
\bigcup_{j=1}^k Y_j \cup \bigcup_{j=1}^k Z_j$, and $K$ intersects all these
intervals. Let $f:\RR\to\RR$ be a contraction of ratio at most $r$. Denote by
$m_Y$ (resp.~$m_Z$) the number of intervals $Y_j$ (resp.~$Z_j$) 
met by some $f(I)$, where
$I$ ranges over the $X_i$'s, $Y_j$'s and $Z_j$'s. We will be done once we show
that $m_Y<k$ or $m_Z<k$.

Using
$\e_1\le\frac{d}{4}$ and (\ref{eq:dist}) we obtain
\beq\lab{eq:kd}
diam\left(f(\bigcup_{j=1}^k Y_j)\right) < diam\left(\bigcup_{j=1}^k
Y_j\right) = (k-1)d+2\e_1 \le kd-2\e_1 \le 
\eeq
\[
\le dist(\{y_1,\dots,y_k\},\{z_1,\dots,z_k\}) - 2\e_1 =
dist(\bigcup_{j=1}^k Y_j,\bigcup_{j=1}^k Z_j), 
\]
so $f\left(\bigcup_{j=1}^k Y_j\right)$ cannot intersect both $\bigcup_{j=1}^k
Y_j$ and $\bigcup_{j=1}^k Z_j$. 
Of course, the same holds for $f\left(\bigcup_{j=1}^k Z_j\right)$, so without
loss of generality we may assume that 
\beq\lab{eq:disj}
\left(\bigcup_{j=1}^k Y_j\right) \cap f\left(\bigcup_{j=1}^k Y_j\right) =
\eset \textrm{ or } \left(\bigcup_{j=1}^k Y_j\right) \cap f\left(\bigcup_{j=1}^k
Z_j\right) = \eset. 
\eeq
For $1\le j_1 <j_2 \le k$ we have $dist(Y_{j_1}, Y_{j_2}) \ge d-2\e_1 \ge
2\e_1$, so if $I$ is an interval of length $2\e_1$ then $f(I)$ cannot intersect 
both $Y_{j_1}$ and $Y_{j_2}$. Moreover, if $H\su\RR$ intersects $t$ many distinct
$Y_j$ intervals, then clearly $diam(H) > d(t-1)-2\e_1 > d(t-1)-d = d(t-2)$,
hence
\beq\lab{eq:diam}
t<\frac{diam(H)}{d}+2.
\eeq
We would like to apply this to $f\left(\bigcup_{j=1}^k Y_j\right)$ and
$f\left(\bigcup_{j=1}^k Z_j\right)$. Clearly
\[
diam\left(f(\bigcup_{j=1}^k Y_j)\right) \leq r diam\left(\bigcup_{j=1}^k
Y_j\right) = r[(k-1)d+2\e_1] \le rkd, 
\]
so by (\ref{eq:diam}) $f\left(\bigcup_{j=1}^k Y_j\right)$ can only meet at
most $rk+2$ many $Y_j$'s, and similarly for $f\left(\bigcup_{j=1}^k
Z_j\right)$. In fact, by (\ref{eq:disj}) we only need to calculate with one of
these two amounts, and altogether we obtain
\[
m_Y < rk+2+n,
\]
where $n$ comes from the $X_i$'s. But by (\ref{eq:k}) $rk+2+n < k$, which
finishes the proof.
\ep

\br
Note that if every contraction $f:X\to X$ is constant, then $X$ clearly has
the Banach Fixed Point Property.
\er

\bt
There exists a non-closed $G_\de$ set $X\su\RR$ with the Banach Fixed Point
Property. Moreover, $X\su [0,1]$ and every contraction mapping $X$ into itself
is constant. 
\et

\bp
Let $K\su\RR$ be a nonempty compact set such that no portion of $K$ can be
covered by a contractive image of $K$. Then $K$ is clearly nowhere dense. 
Let
\[
X = (K+\QQ)^c \cap [0,1],
\]
then $X$ is $G_\de$. As $K+\QQ$ is a nonempty set of the
first category, it is not open in $[0,1]$, hence $X$ is not closed.

Now, let $f:X\to X$ be a non-constant
contraction. As above, let $\tilde{f}:\RR\to\RR$ be a
contraction extending $f$. As $X$ is dense in $[0,1]$, we have
$\tilde{f}([0,1]) \su [0,1]$. We can clearly assume that $\tilde{f}$ is
constant on $(-\infty, 0]$ and $[1,\infty)$, hence $ran(\tilde{f}) \su
[0,1]$. Then $ran(\tilde{f})$ is a nondegenerate interval 
$I\su [0,1]$. Pick $q_0\in\QQ$ so 
that $(K+q_0) \cap int(I) \neq \eset$. As $\tilde{f}(X)\su X$, we have $X^c
\cap I \su \tilde{f}(X^c)$, so $(K+q_0) \cap I \su
(K+\QQ) \cap I \su \tilde{f}(K+\QQ) \cup \{\tilde{f}(0),\tilde{f}(1)\}=
\bigcup_{q\in\QQ} \tilde{f}(K+q) \cup \{\tilde{f}(0),\tilde{f}(1)\}$. Since
$K$ is nowhere dense, there is a nondegenerate interval
$[a,b]\su int(I)$ intersecting $K+q_0$ such that $a,b\notin K+q_0$. The closed
set $[a,b] \cap (K+q_0) \su \bigcup_{q\in\QQ} \tilde{f}(K+q) \cup
\{\tilde{f}(0),\tilde{f}(1)\}$, which is a 
covering by countably many closed sets, hence by the Baire Category Theorem
one of them covers a portion of $K+q_0$, which contradicts the choice of $K$.
\ep

\bt
There exists a non-closed $F_\si$ subset of $[0,1]$ with the Banach Fixed Point
Property.
\et

\bp
Again, let $K\su\RR$ be a nonempty nowhere dense compact set such that no
portion of $K$ can be covered by a contractive image of $K$. Then clearly $K$
has no isolated points, so $K$ is homeomorphic to the Cantor set
\cite[7.4]{Ke}. We can clearly assume that $min(K)=0$ and $max(K)=1$.
Let $\{I_n\}_{n\in\NN}$ be the set of contiguous open
intervals of $K$. Set
\[
X = \cup_{n\in\NN} \overline{I_n}.
\]
That is, $X$ is `$[0,1] \sm K$ plus the endpoints'. 
This set is clearly $F_\si$, and it is not closed, as it is dense in $[0,1]$
but only contains countably many points of $K$.   

In order to show that it has the BFPP let $f:X\to X$ be a
contraction, and as above, let $\tilde{f}:\RR\to [0,1]$ be a
contraction extending $f$ (here we use again that $X$ is dense in $[0,1]$) that
is constant on $(-\infty, 0]$ and $[1,\infty)$. 
If $\tilde{f}$ is constant then we are done,
otherwise $ran(\tilde{f})$ is a nondegenerate interval $I\su [0,1]$. 
If $I\su X$ then 
(by connectedness) we have $I\su \overline{I_{n_0}}$ for some $n_0\in\NN$, and
therefore $f\rest_{\overline{I_{n_0}}}$ has a fixed point. 

So we can assume $X^c\cap I \neq \eset$. Then using again that $X$ is a union
of closed intervals we obtain that $X^c\cap int(I) \neq
\eset$. Choose a nondegenerate interval $[a,b]\su int(I)$ intersecting $K$ so
that $a,b \notin K$. Similarly as
above, $X^c\cap I \su \tilde{f}(X^c) \su \tilde{f}(K)$. As this last set is
closed, $\overline{X^c\cap I} \su \tilde{f}(K)$. Set $E=\bigcup_{n\in\NN}
\left( \overline{I_n}\sm I_n \right)$; that is, the set of
endpoints. Then $K \cap [a,b] = \overline{(K \sm E)\cap [a,b]} =
\overline{X^c\cap [a,b]} \su \overline{X^c\cap I} \su \tilde{f}(K)$, which
is impossible by the choice of $K$.
\ep

It is well known \cite[3.11]{Ke} that there is a complete metric
equivalent to the usual one on a set $X\su\RR^n$ iff $X$ is $G_\de$.
Combining this fact with the above theorem and Theorem \ref{amb} 
we obtain the following.

\bcor
There is a bounded Borel (even $F_\si$) subset of $\RR$ with the Banach Fixed
Point Property that is not complete with respect to any equivalent metric.
\ecor

Finally we show that even a nonmeasurable set can have the BFPP. A set
$B\su [0,1]^n$
is called a \emph{Bernstein set} if $B\cap F\neq \eset$ and $B^c\cap F\neq
\eset$ for every uncountable closed set $F\su[0,1]^n$. 
It is well known that every Bernstein set is nonmeasurable \cite[5.3]{Ox}
(which works for $[0,1]^n$ instead of $\RR$).

\bt
For every integer $n>0$ there exists a nonmeasurable set in $\RR^n$ 
with the Banach Fixed Point Property. Moreover, there exists a
Bernstein set in $[0,1]^n$ with the BFPP, such that every contraction mapping
this set into itself is constant.
\et

\bp
It suffices to prove the second statement.
Enumerate the uncountable closed sets $F\su [0,1]^n$ as $\{F_\al:\al<2^\om\}$,
and also the non-constant contractions $f:\RR^n\to\RR^n$ as
$\{f_\al:\al<2^\om\}$. We define a characteristic function
$\ph:[0,1]^n\to\{0,1\}$, and the Bernstein set with the required properties will
be $X = \{x\in[0,1]^n:\ph(x)=1\}$.

Suppose we
have already defined $\ph$ on a set $D_\al\su [0,1]^n$ of cardinality
$<2^\om$. 
We define it for four more points.
As every uncountable closed set is of cardinality $2^\om$, we can
pick two distinct points $x_\al, y_\al\in F\sm D_\al$ and define 
$\ph(x_\al) = 0$, $\ph(y_\al) = 1$. This will make sure that $X$ will be a
Bernstein set in $[0,1]^n$.

As $ran(f_\al)$ is a nondegenerate connected set, its projection on every line
is an interval, and for a suitable line this interval is nondegenerate. Hence
$|ran(f_\al)|=2^\om$. Therefore $|ran(f_\al)\sm
(D_\al\cup\{x_\al,y_\al,Fix(f_\al)\})| = 2^\om$, where $Fix(f_\al)$ is the
(unique) fixed point of $f_\al$. 
As the inverse images of the points of this
set form a disjoint family of size $2^\om$ of nonempty sets, and
$|D_\al\cup\{x_\al,y_\al,Fix(f_\al)\}| < 2^\om$, 
there exists $u_\al \in ran(f_\al)\sm (D_\al \cup \{x_\al,y_\al,Fix(f_\al)\})$
such that $f_\al^{-1}(u_\al) \cap (D_\al\cup\{x_\al,y_\al,Fix(f_\al)\}) =
\eset$. Pick an arbitrary $v_\al \in f_\al^{-1}(u_\al)$, then $v_\al \neq
u_\al$. Finally, define $\ph(u_\al) = 0$, $\ph(v_\al) = 1$.

After finishing the induction define $\ph$ to be $0$ outside
$\bigcup_{\al<2^\om} D_\al$. As we mentioned above, $X$ is easily seen to be a
Bernstein set in $[0,1]^n$. In order to get a contradiction, let $f:X\to X$ be a
non-constant contraction. Then it can be extended to $\RR^n$, so $f=f_\al$ for
some $\al$. But then $v_\al\in X$ and $f(v_\al) = f_\al(v_\al) = u_\al \notin
X$, a contradiction.
\ep

\br
It is not hard to see that if $X = sin(1/x)\rest_{(0,1]}$ then there exists a
function $f:X\to X$ with no fixed points such that $|f(x)-f(y)|<|x-y|$ for
every $x,y\in X$. (Just `map each wave horizontally to the next one'.)
It would be interesting to know what happens if we replace the class of
contractions with this larger class of strictly distance-decreasing functions.
\er

\bq
Is there for some $n\in\NN$ a non-closed $F_\si$ subset $X\su\RR^n$ with the
Banach Fixed Point Property such that every contraction $f:X\to X$ is
constant? Is there such a simultaneously $F_\si$ and $G_\de$ set?
\eq

\noindent
\textbf{Acknowledgement} The author is indebted to T.~Keleti and M.~Laczkovich
for some helpful discussions.

\bigskip

\noindent
\textsc{M\'arton Elekes} 

\noindent
\textsc{R\'enyi Alfr\'ed Institute of Mathematics} 

\noindent
\textsc{Hungarian Academy of Sciences} 

\noindent
\textsc{P.O. Box 127, H-1364 Budapest, Hungary}

\noindent
\textit{Email:} \verb+emarci@renyi.hu+

\noindent
\textit{URL:} \verb+www.renyi.hu/~emarci+


\begin{thebibliography}{aaaa}

\bibitem{Ba} A.~C.~Babu, 
A converse to a generalised Banach contraction principle,
\textsl{Publ.~Inst.~Math.~(Beograd) (N.S.)} \textbf{32(46)}, (1982), 5--6.

\bibitem{Be1} E.~Behrends, Problem Session of the 34th Winter School in
  Abstract Analysis, 2006.

\bibitem{Be2} E.~Behrends, private communication, 2006.

\bibitem{Bes} C.~Bessaga, On the converse of the Banach "fixed-point
principle", \textsl{Colloq.~Math.} \textbf{7}, (1959), 41--43. 

\bibitem{Ci} L.~B.~\'Ciri\'c, On some mappings in metric spaces and fixed
points, \textsl{Acad.~Roy.~Belg.~Bull.~Cl.~Sci.~(6)} \textbf{6}, (1995),
no.~1-6, 81--89.

\bibitem{Fe} H.~Federer: \textsl{Geometric Measure Theory.} Springer-Verlag,
1969.

\bibitem{Iv} A.~A.~Ivanov,
Fixed points of mappings of metric spaces, \textsl{Studies in topology,
II.~Zap.~Nau\v cn.~Sem.~Leningrad.~Otdel.~Mat.~Inst.~Steklov.~(LOMI)}
\textbf{66}, (1976), 5--102, 207.

\bibitem{Ja} J.~Jachymski, General solutions of two functional inequalities
and converses to contraction theorems, \textsl{Bull.~Polish Acad.~Sci.~Math.}
\textbf{51} (2003), no.~2, 147--156.

\bibitem{Ja1} L.~Jano\v s, A converse of Banach's contraction theorem,
  \textsl{Proc. Amer. Math. Soc.}, \textbf{18}, (1967), 287--289.

\bibitem{Ja2} L.~Jano\v s, 
A converse of the generalised Banach's contraction theorem,
\textsl{Arch.~Math.~(Basel)} \textbf{21}, (1970), 69--71.

\bibitem{Ke} A.~S.~Kechris: \textsl{Classical Descriptive Set Theory.}
Springer-Verlag, 1995. 

\bibitem{Ki} W.~A.~Kirk, Contraction mappings and extensions, \textsl{Handbook
of metric fixed point theory,} 1--34, Kluwer Acad.~Publ., Dordrecht, 2001. 

\bibitem{Ku} K.~Kuratowski: \textsl{Topology.} Academic Press, 1966.

\bibitem{Me} P.~R.~Meyers, A converse to Banach's contraction theorem,
\textsl{J. Res. Nat. Bur. Standards Sect. B} \textbf{71B}, (1967), 73--76.

\bibitem{MP} A.~Mukherjea and K. Pothoven:
\textsl{Real and functional analysis.
Mathematical Concepts and Methods in Science and Engineering, Vol. 6.}
Plenum Press, New York-London, 1978.

\bibitem{Op} V.~I. Opo\u\i cev, 
A converse of the contraction mapping principle, \textsl{Uspehi Mat. Nauk}
\textbf{31}, (1976), no. 4 (190), 169--198.

\bibitem{Ox} J.~C.~Oxtoby: \textsl{Measure and Category. A survey of the
analogies between topological and measure spaces.} Second edition.
Graduate Texts in Mathematics No. 2, Springer-Verlag, 1980.

\bibitem{Ru} I.~A.~Rus: \textsl{Generalised contractions and applications.}
Cluj University Press, Cluj-Napoca, 2001.

\end{thebibliography}
\end{document}